\documentclass{amsart}
\usepackage{amscd,amssymb}

\theoremstyle{plain}
\newtheorem{thm}{Theorem}[section]
\newtheorem{lem}[thm]{Lemma}
\newtheorem{definicja}[thm]{Definition}

\newtheorem{prop}[thm]{Proposition}

\theoremstyle{definition}

\theoremstyle{remark}

\numberwithin{equation}{section}

\begin{document}
%\setcounter{section}{-1}
%Topmatter
\title{Polish group actions and computability} 
\thanks{{\bf 2000 Mathematics Subject Classsification}: 03E15, 03D45}
\author{ Aleksander Ivanov and Barbara Majcher-Iwanow}

\maketitle

\bigskip

\bigskip

\begin{quote}
ABSTRACT. 
Let $G$ be a closed subgroup of $S_{\infty}$ and ${\bf X}$ be a Polish 
$G$-space with a countable basis $\mathcal{A}$ of clopen sets. 
Each $x\in {\bf X}$ defines a characteristic function $\tau_x$ 
on $\mathcal{A}$ by $\tau_x (A)=1 \Leftrightarrow x\in A$.  
We consider computable complexity of $\tau_x$ and some related questions. 
\end{quote}
\bigskip

\section{Introduction}

Let $L= (R^{n_i}_i)_{i\in I}$ be a countable relational language
and $\mathbf{X}_L=\prod_{i\in I}2^{\omega^{n_i}}$ be the corresponding
topological space under the product topology.
We consider $\mathbf{X}_L$ as the space of all $L$-structures on $\omega$
(see Section 2.5 in \cite{bk} or Section 2.D of \cite{becker}
for details).
If $F$ is a countable {\em fragment}
\footnote{here we assume that $F$ is closed under $\wedge$, $\vee$ and $\neg$, 
and do not assume that $F$ is closed under quantifiers or subformulas} 
of $L_{\omega_1 \omega}$, then the family of all sets
$Mod(\phi ,\bar{s})=\{ M\in \mathbf{X}_L :M\models \phi (\bar{s})\}$,
where $\phi\in F$ and $\bar{s}$ is a tuple from $\omega$, forms a basis of
a topology on $\mathbf{X}_L$ which will be denoted by $\mathbf{t}_{F}$ 
(it is easy to see that the fragment of quantifier-free 
first-order formulas defines the original product topology). 
The group $S_{\infty}$ of all permutations of $\omega$ has 
the natural action on $\mathbf{X}_L$ and the action is continuous 
with respect to $\mathbf{t}_F$. 
It is called the {\em logic action} of $S_{\infty}$ 
on $(\mathbf{X}_L ,\mathbf{t}_F )$. 
For $M\in \mathbf{X}_L$ we define the characteristic function 
$\tau_M$ distinguishing in the above basis of the topology $\mathbf{t}_F$, 
clopen sets containing $M$. 
Using the standard coding of terms and formulas we see that computable 
complexity of $\tau_M$ corresponds to complexity of $M$ studied 
in computability theory.  
{\em The aim of our paper is to show that this idea extends the approach 
of computability theory to Polish group actions and nice topologies} 
(introduced in \cite{becker2}). 
In particular we show that decidable theories can be considered  
as so called {\em decidable pieces of the canonical partition}. 
Identifying such pieces with appropriate computable functions 
we now consider complexity of some natural properties of pieces, 
for example  counterparts of $\omega$-categoricity. 
In particular we develope and generalize some material from \cite{LS} 
concerning complexity of the family of $\omega$-categorical theories. 
\parskip0pt        

We illustrate our approach by some adaptations of examples 
of non-G-compact theories from \cite{sasza} and \cite{P}. 
We have found that they also provide some new theories 
having (having no) degree. 
This material is also given in the general form of 
Polish $G$-spaces. 
The final part of the paper contains new examples of groups 
with and without degrees. 
These groups are $\omega$-categorical. 

To present our approach in more detail we should 
remind the reader some definitions. 
In particular we must explain what a nice topology is.

\subsection{Preliminaries.} 

A {\em Polish space (group)} is a separable, completely
metrizable topological space (group).
If a Polish group $G$ continuously acts on a Polish space ${\bf X}$,
then we say that ${\bf X}$ is a {\em Polish $G$-space}.
We usually assume that $G$ is considered under a left-invariant metric.
We simply say that a subset of ${\bf X}$ is {\em invariant} 
if it is $G$-invariant. 
\parskip0pt

We consider the group $S_{\infty}$ of all permutations
of the set $\omega$ of natural numbers under the usual
left invariant metric $d$ defined by
$$
d(f,g)=2^{-\min\{k:f(k)\not=g(k)\}},\mbox{ whenever } f\not=g.
$$
For a finie set $D$ of natural numbers let $id_D$ be the identity map 
$D\rightarrow D$ and $V_{D}$ be the group of all
permutations stabilizing $D$ pointwise, i.e.,
$V_D =\{f\in S_{\infty}: f(k)=k\mbox { for every } k\in D\}$. 
Writing $id_{n}$ or $V_n$  we treat $n$ as the set of 
all natural numbers less than $n$. \parskip0pt 

Let $S_{<\infty}$ denote the set of all bijections between 
finite substes of $\omega$. 
We shall use small greek letters $\delta, \sigma, \tau$ 
to denote elements of $S_{<\infty}$. 
For any  $\sigma \in S_{<\infty}$ let
$dom[\sigma], rng[\sigma]$ denote the domain
and the range of $\sigma$ respectively.

For every $\sigma \in S_{<\infty}$ let
$V_{\sigma} =\{f\in S_{\infty}:f\supseteq \sigma \}$.
Then for any $f\in  V_{\sigma}$
we have  $V_{\sigma}=fV_{dom[\sigma]}=V_{rng[\sigma]}f$.
Thus the family $\mathcal{N}=\{ V_{\sigma}:\sigma\in S_{<\infty}\}$ 
consists of all left (right) cosets of all subgroups $V_D$ as above. 
This is a basis of the topology of $S_{\infty}$. 

Given $\sigma\in S_{<\infty}$ and $D\subseteq dom[\sigma]$ 
for any $f\in V_{\sigma}$ we have $V_D^f=V_{\sigma [D]}$,
where $V_D^f$ denotes the conjugate $fV_D f^{-1}$.
\parskip0pt

In our paper we concentrate on Polish $G$-spaces, 
where $G$ is a closed subgroup of $S_{\infty}$.  
For such a group we shall use 
the relativized version of the above, i.e.,
$V^G_{\sigma}=\{ f\in G:f\supseteq \sigma\}$,
$S_{<\infty}^G=\{ f|_D : f\in G$ and 
$D$ is a finite set of natural numbers $\}$ 
(observe that for any subgroup $G$ and any finite set $D$ 
of natural numbers we have $id_{D}\in S_{<\infty}^G$)
and $V^G_{\sigma}= V_{\sigma}\cap G$. 
The family 
$\mathcal{N}^G =\{ V_{\sigma}^G: \sigma\in S_{<\infty}^G\}$ 
is a basis of the standard topology of $G$. \parskip0pt 

All basic facts concerning Polish $G$-spaces can 
be found in \cite{bk}, \cite{hjorth} and \cite{kechris}. 
\parskip0pt 

Since we will use Vaught transforms, 
recall the corresponding definitions. 
The Vaught $*$-transform of a set $B\subseteq {\bf X}$ 
with respect to an open $H\subseteq G$ is the set
$B^{*H}=\{ x\in X:\{ g\in H:gx\in B\}$ is comeagre in $H\}$.
We will also use another Vaught transform
$B^{\Delta H}=\{ x\in X:\{ g\in H:gx\in B\}$ 
is not meagre in $H\}$. 
It is worth noting that
{\em for any open $B\subseteq X$ and any open $K<G$ we have}
$B^{\Delta K}=KB$. 
Indeed, by continuity of the action for any $x\in KB$ and
$g\in K$ with $gx\in B$ there are open neighbourhoods
$K_1\subseteq K$ and $B_1\subseteq KB$ of $g$ and $x$
respectively so that
$K_1 B_1\subseteq B$; thus $x\in B^{\Delta K}$.
Other basic properties of Vaught transforms can be found 
in \cite{bk}. 

\subsection{Nice bases.} 

We now define nice topologies.  
Let $G$ be a closed subgroup of $S_{\infty}$ and 
let $(\langle {\bf X}, \tau\rangle ,G)$ be 
a Polish $G$-space with a countable basis $\mathcal{A}$.
Along with the topology $\tau$ we shall consider 
another topology on ${\bf X}$.
The following definition comes from \cite{becker2}.

\begin{definicja} \label{nto}
A topology $\mathbf{t}$ on ${\bf X}$ is {\em nice} for the $G$-space
$(\langle {\bf X}, \tau\rangle , G)$ if the following
conditions are satisfied.\\
(a) $\mathbf{t}$ is a Polish topology, $\mathbf{t}$ is finer than $\tau$ 
and the $G$-action remains continuous with respect to $\mathbf{t}$.\\
(b) There exists a basis $\mathcal{B}$ for $\mathbf{t}$ such that:\parskip0pt

(i) $\mathcal{B}$ is countable;\parskip0pt

(ii) for all $B_1, B_2 \in \mathcal{B}$, $B_1\cap B_2\in \mathcal{B}$;
\parskip0pt

(iii) for all $B\in \mathcal{B}$, ${\bf X}\setminus B\in \mathcal{B}$; 
\parskip0pt

(iv) for all $B\in \mathcal{B}$ and $u\in \mathcal{N}^G$,
$B^{*u} \in \mathcal{B}$;\parskip0pt

(v) for any $B\in \mathcal{B}$ there exists an open subgroup $H < G$
such that $B$ is invariant \parskip0pt 

under the corresponding $H$-action.\\
A basis satisfying condition $(b)$ is called a {\em nice} basis.
\end{definicja}

In this definition $B^{*u}$ denotes the Vaught
$*$-transform of $B$. 
It is noticed in \cite{becker2} that any nice basis also
satisfies property (b)(iv) of the definition above
for $\Delta$-transforms. 
As we have already mentioned above, 
for any $B\in \mathcal{B}$ and any open $K<G$ we have 
$B^{\Delta K}=K\cdot B$. \parskip0pt

From now on $\mathbf{t}$ will always stand for a nice topology 
on ${\bf X}$ and $\mathcal{B}$ will be its nice basis.
Observe that
{\em any nice basis is invariant in the sense that
for every $g\in G$ and $B\in \mathcal{B}$ we have
$gB\in \mathcal{B}$}.
Indeed, by (v), there is $u\in \mathcal{N}^G$
such that $B$ is $u$-invariant.
Using  properties of Vaught transforms, we obtain
the equalities $gB=gB^{*u}=B^{*ug^{-1}}$.
Then we are done by (iv). \parskip0pt 

By Theorem 1.11 from \cite{becker2} for any $G$-space 
$({\bf X},\tau )$ as in Definition \ref{nto} a nice topology 
$\mathbf{t}$ always exists.     
In our paper  we will be interested in nice topologies ${\bf t}$ 
such that $\mathcal{B}_{\mathbf{t}}$ 
is effectively coded.  \parskip0pt 

Nice bases naturally arise when we consider the situation 
described in the beginning of our introduction. 
Let $L$ be a countable relational language
and $\mathbf{X}_L$ be the corresponding $S_{\infty}$-space under 
the product topology $\tau$ and the corresponding logic 
action of $S_{\infty}$. 
Let $\mathbf{t}_F$ be the topology on $\mathbf{X}_L$ corresponding to 
some countable fragment of $L_{\omega_1 \omega}$-formulas 
as it was described above. 
Theorem 1.10 of \cite{becker2} states that if $F$ is closed 
with respect to quantifiers, then $\mathbf{t}_F$ is nice.  
In this case usually the basis defining 
$\mathbf{t}_F$ is effectively coded.

\section{Polish group actions and decidable relations}

\subsection{Approach}

Our circumstances are standard and in particular, arise 
when one studies $S_{\infty}$-spaces of logic actions. 
Let $G$ be a closed subgroup of $S_{\infty}$ and 
$({\bf X},\tau )$ be a Polish $G$-space. 
Let $\mathcal{A}$ be a countable basis of $({\bf X},\tau )$  
closed with respect to $\cap$. 
We assume that each $A$ of $\mathcal{A}$ is $H$-invariant with
respect to some basic subgroup $H\in \mathcal{V}^G$. 
We will also assume that 
the subfamily of $\mathcal{A}$ 
consisting of clopen sets generates the same topology. 
\parskip0pt 

We assume that the bases $\mathcal{N}^G$ and $\mathcal{A}$ 
are computably 1-1-enumerated so that the relations of 
inclusion $\subseteq$ together with the corresponding 
operations $\cap$ (as well as the predicates $Clopen$ 
for the set of clopen subsets of $\mathcal{A}$ and 
$\mathcal{V}^G$ for the set of all basic subgroups 
from $\mathcal{N}^G$ respectively) are presented 
by decidable relations on $\omega$. 
Moreover we assume that there is an algorithm deciding 
the problem if for a basic clopen set $U$ 
(of $\mathcal{N}^G$ or $\mathcal{A}$) and a natural 
number $i$ the diametr of $U$ is less than $2^{-i}$.  

We also assume that the following relations are decidable: \parskip0pt 

(a) 
$Inv(V,U)\Leftrightarrow (V\in \mathcal{V}^G)\wedge (U\in \mathcal{A})\wedge (U$ 
is $V$-invariant $)$ ; \parskip0pt 

(b) 
$Orb_{m,n}(N,V_1 ,...,V_m ,V_{m+1},...,V_{2m},U_1 ,...,U_n ,U_{n+1},...,U_{2n})\Leftrightarrow$ 
$(N \in \mathcal{N}^G )\wedge\bigwedge^{2m}_{i=1}(V_i \in \mathcal{V}^G )\wedge$
$\bigwedge^{2n}_{i=1}(U_i \in\mathcal{A})\wedge ($ 
the tuple $(V_{m+1},...,V_{2m},U_{n+1},...,U_{2n})$ is of
the form $(V^{g}_{1},...,V^{g}_m ,gU_1 ,...,gU_n )$  
for some $g\in N)$.

\begin{definicja} \label{Comp} 
We say that an element $x\in {\bf X}$ is {\em computable} if the relation 
$$
Sat_x (U)\Leftrightarrow (U\in \mathcal{A})\wedge (x\in U)
$$ 
is decidable. 
\end{definicja} 

In the case of the logic action, when $x$ is a structure 
on $\omega$, this notion is obviously equivalent to 
the notion of a computable structure.  
We will denote by $Sat_x ({\mathcal{A}})$ the set 
$\{ C\in {\mathcal{A}}: Sat_x (C)$ holds $\}$. 
It is straightforward that 
\begin{quote}
for a computable $x$ there is a computable function 
$\kappa :\omega \rightarrow \mathcal{A}$ such that for all 
natural numbers $n$, $x\in \kappa (n)$ and $\kappa (n)$ is 
clopen with $diam(\kappa (n))\le 2^{-n}$. 
\end{quote} 
It is also worth noting that when $\mathcal{A}$ consists 
of clopen sets, the existence of such a computable function 
$\kappa$ already implies that the relation $Sat_x$ is decidable. 
Indeed, since $A$ is clopen, in order to decide $Sat_x (A)$ 
we have to check if  $(\exists l)(\kappa (l)\subset A)$ 
or $(\exists l)(\kappa (l)\cap A=\emptyset )$. 
 
We also say that an element $g\in G$ is {\em computable} 
if the relation $(N\in \mathcal{N}^G)\wedge (g\in N)$ 
is computable. 
Then there is a computable function realizing the 
same property as $\kappa$ above but already in the case 
of the basis $\mathcal{N}^G$. 
Since $\mathcal{N}^G$ consists of clopen sets 
these two properties are equivalent. 
In the following lemma we use standard indexations 
of the set of computable functions and of the set 
of all finite subsets of $\omega$. 

\begin{lem} \label{21}
The following relations belong to $\Pi^0_2$:\\
(1) $\{ e:$ the function $\varphi_e$ is a characteristic 
function of a subset of $\mathcal{A}\}$; \\ 
(2) $\{ (e,e'):$ there is a computable element $x\in {\bf X}$ 
such that the function $\varphi_e$ is a characteristic function 
of the set $Sat_x ({\mathcal{A}})$ and the function $\varphi_{e'}$ 
realizes the corresponding function $\kappa$ defined 
after Definition \ref{Comp} $\}$;\\   
(3) $\{ (e,e'):$ there is an element $g\in G$ such that the 
function $\varphi_e$ is a characteristic function of the 
subset $\{ N\in\mathcal{N}^G:g\in N\}$ and the function 
$\varphi_{e'}$ realizes the corresponding function $\kappa$ 
defined after Definition \ref{Comp} (in the case of $\mathcal{N}^G$) $\}$.    
\end{lem} 

{\em Proof.} (1) Obviuos. 
Here and below we use the fact that a function is computable if 
and only if its graph is computably enumerable. 

(2) The corresponding definition can be described as follows:
 
$$
("e \mbox{ is a characteristic function of a subset of }\mathcal{A}")\wedge 
$$ 
$$
(\forall n)(Clopen(\varphi_{e'}(n))\wedge (\varphi_{e'}(n)\not=\emptyset )
\wedge (\varphi_{e}(\varphi_{e'}(n))=1 )\wedge 
diam(\varphi_{e'}(n))<2^{-n} ) \wedge  
$$
\begin{quote} 
$(\forall d)(\exists n )(($ 
"every element $U'$ of the finite subset of $\mathcal{A}$ 
with the canonical index $d$ satisfies $\varphi_e (U') =1$") 
$\leftrightarrow ($ "$\varphi_{e'}(n)$ is contained in any element $U'$ 
of the finite subset of $\mathcal{A}$ with the canonical index $d$"$))$.  
\end{quote}
It is clear that by Cantor's theorem the last part of the 
conjunction ensures the existence of the corresponding $x$. 

(3) is similar to (2). 
$\Box$
\bigskip 

We now describe how decidability of elementary 
theories appears in our approach.  
 
By Proposition 2.C.2 of \cite{becker} there exists a unique
partition of ${\bf X}$, ${\bf X}=\bigcup\{ Y_{t}: t\in T\}$, 
into invariant $G_{\delta}$-sets $Y_{t}$ such that every 
$G$-orbit from $Y_{t}$ is dense in $Y_{t}$.
It is called the {\em canonical partition} of the $G$-space ${\bf X}$.
To construct this partition take $\{ A_{j}\}$, a countable basis
of ${\bf X}$, and for any $t\in 2^{\omega}$ define
\[
Y_{t}=(\bigcap\{ GA_{j}:t(j)=1\})\cap
(\bigcap\{ {\bf X}\setminus GA_{j}:t(j)=0\})
\]
and take $T=\{ t\in 2^{\omega}:Y_{t}\not=\emptyset\}$.

We say that a piece $Y_{t}$ is {\em decidable} if 
the corresponding function $\mu_t :\omega\rightarrow 2$ 
characterizing all $A_j$ with $Y_t \subseteq GA_j$, is computable. 

In the case of the logic action of $S_{\infty }$ on the space
${\bf X}_{L}$ of countable $L$-structures under the topology 
$\mathbf{t}_F$ (corresponding to a fragment $F$; see Introduction),
each piece of the canonical partition is an equivalence class
with respect to the $F$-elementary equivalence 
$\equiv_{F}$ \cite{becker}. 
Thus a computable piece is a decidable complete $F$-elementary theory.  

We apply this idea to nice topologies corresponding to $({\bf X},\tau )$. 

\begin{definicja} 
Let $\mathcal{B}$ be a nice basis corresponding to a nice topology 
$\mathbf{t}$ of $({\bf X},\tau )$. 
We say that the basis $\mathcal{B}$ is {\em computable} if 
$\mathcal{B}$ is computably 1-1-enumerated so that 
there is a computable function $\mathcal{A}\rightarrow \mathcal{B}$ 
finding the $\mathcal{B}$-numbers of elements of $\mathcal{A}$ 
(such that $\mathcal{A}$ is computable) 
and the following relations are decidable: 

(i) the binary relations of inclusion $\subseteq$, and taking 
the complement: $B'={\bf X}\setminus B$; 

(ii) binary relation 
$Inv(V,U)\Leftrightarrow (V\in \mathcal{V}^G)\wedge (U\in \mathcal{B})\wedge (U$ 
is $V$-invariant $)$;

(iii) ternary relations corresponding to the operation $\cap$ 
($B_1 \cap B_2 =B_3$) and the operation of taking the Vaught 
transforms : $B^{*u}_1 =B_2$ and $B^{\Delta u}_1 =B_2$. 
\end{definicja} 
 
Using the same definition as above we can define decidable 
pieces of the canonical partition corresponding to $\mathcal{B}$. 
On the other hand since for every $A\in \mathcal{A}$ the 
element $GA=A^{\Delta G}$ belongs to $\mathcal{B}$, each $\tau$-canonical 
piece is an intersection of an appropriate subset of $\mathcal{B}$.  
Now $\tau$-canonical pieces become more tractable.  

\begin{prop} \label{22} 
Let $\mathcal{B}$ be a computable nice basis corresponding 
to a nice topology $\mathbf{t}$ of $({\bf X},\tau )$. \\
(1)  The following relation belongs to $\Pi^0_2$:  

$\{ (e,e',e'',A): A\in \mathcal{A}$ and there is 
a computable element $x\in A$ such that 

the function $\varphi_e$ is a characteristic 
function of the set $Sat_x ({\mathcal{A}})$, 

the function $\varphi_{e'}$ realizes the corresponding 
function $\kappa$ as after Definition \ref{Comp}, 
 
$\varphi_{e''}$ is a characteristic function on 
$\mathcal{A}$ defining a piece of the canonical partition, 

and the computable element $x$ belongs to the canonical 
piece defined by $\varphi_{e''} \}$. \\ 
(2) The class $\Pi^0_4$ contains the set of all $e''$ such that 
$\varphi_{e''}$ codes a decidable piece of the $\tau$-canonical 
partition such that all computable elements of the piece are 
contained in the same orbit of computable elements of $G$.  
\end{prop} 

{\em Proof.} 
(1) By Lemma \ref{21}(2) the statement that $\varphi_{e}$ 
and $\varphi_{e'}$ realize a computable element $x$ from 
${\bf X}$, belongs to $\Pi^0_{2}$. 
As in Lemma \ref{21}(1) we see that the statement that 
$\varphi_{e''}$ is a characteristic function 
on $\mathcal{A}$, also belongs to $\Pi^0_2$. 
Since $\tau$ is generated by clopen members of $\mathcal{A}$, 
to express that $x$ belongs to the intersection of $A$ and the canonical 
piece defined by $\varphi_{e''}$ it suffices to state:  
\begin{quote} 
(a) $(\exists l)(\varphi_{e'}(l)\subseteq A)$, \\  
(b) for any $l$ and elements $B_1 ,...,B_k$ of 
$\mathcal{A}$ the intersection  
$$
\bigcap \{ GB_i :\varphi_{e''}(B_i )=1\} \cap\varphi_{e'}(l) \cap 
\bigcap  \{ {\bf X}\setminus GB_i :\varphi_{e''}(B_i )=0\} 
$$ 
is non-empty and \\ 
(c) $(\forall B\in \mathcal{A})($ "$\varphi_{e''}(B)=1$ 
is equivalent to  $(\exists l)(\varphi_{e'}(l)\subset GB)$"). 
\end{quote} 
As in the proof of Lemma \ref{21} it is easy to verify 
that these conditions belong to $\Pi^0_2$. 
We also use that $\mathcal{B}$ is a nice basis and the fact 
that $GB=B^{\Delta G}$.  

(2) We express the property of (2) as the statement that 
for any two pairs $(e_1 ,e'_1 )$, $(e_2 ,e'_2 )$ the 
following alternative holds: either one of the tuples 
$(e_1, e'_1 ,e'',{\bf X})$ or $(e_2, e'_2 ,e'',{\bf X})$ 
does not satisfy the condition from (1) or there is 
a number $e_0$ such that $\varphi_{e_0}$ maps $\omega$ 
to a decreasing sequence from $\mathcal{N}^G$ such that for all 
$l,k,k'$  we have $diam(\varphi_{e_0}(l))<2^{-l}$ and  
$\varphi_{e'_1}(k')\cap \varphi_{e_0}(l)\varphi_{e'_2}(k)\not=\emptyset$.    
This is a $\Pi^0_4$-condition. 
$\Box$

\subsection{Compact topologies and G-orbits which are canonical pieces} 

In fact Proposition \ref{22} concentrates on 
"effective parts" of pieces of the canonical partition.  
In this section we make an easy general observation 
(without any neglect of non-computable elements) 
concerning complexity of pieces of the canonical 
partition under the assumption that the basic 
topology $\tau$ is compact. 
The motivation for this assumption is the paper \cite{LS}, 
where it is shown that the complexity of $\omega$-categorical 
first-order theores is $\Pi^0_3$.   
So we concentrate on pieces which are $G$-orbits. 
Following the tradition of computable model theory we 
will restrict ourselves by computable pieces of 
the canonical partition. 
Then each piece can be identified with the corresponding 
computable function (see the previous section). 
Since we do not have some natural logical tools, we 
cannot preserve the statement of \cite{LS} in our context. 
On the other hand we will show that under some natural 
assumptions the level of complexity is very close 
to that of \cite{LS}.  

We start with the following observation. 
\begin{quote}
Let ${\bf t}$ be a nice topology with respect 
to $({\bf X}, \tau, G)$ and $X_0$ be 
a $\tau$-canonical piece. 
If $X_0$ is a $G$-orbit of some $x\in X_0$, then 
both topologies $\tau$ and ${\bf t}$ are equal on $X_0$   
(Proposition 1.4 of \cite{basia}). 
\end{quote} 
On the other hand Theorem 3.4 from \cite{basia} 
(which is a version of Ryll-Nardzewski's 
theorem) states that a ${\bf t}$-canonical piece 
$Y$ is a $G$-orbit if and only if for any basic 
clopen $H<G$ any $H$-type of $Y$ is principal 
(the corresponding terms are defined in \cite{basia}). 
Then a standard logic argument shows that when $X_0$ 
is as above and the induced space $(X_0 ,\tau )$ 
is compact, for any $H\in \mathcal{V}^G$ 
the set of all intersections of $X_0$ with 
$H$-invariant members of 
the nice basis ${\mathcal{B}}$ is finite.  
This allows us to find some counterpart 
of the result from \cite{LS} mentioned above. 
To formulate it we need the following relation. 

We say that $e\in \omega$ and $B\in {\mathcal{B}}$ 
satisfy the relation $Con$ (i.e. $\models Con(e,B)$), 
if there is a decidable $\tau$-canonical piece $Y$ 
such that $B\cap Y\not=\emptyset$, and $\varphi_e$ 
is the characteristic function of the set of all 
$A_j \in {\mathcal{A}}$ with $Y\subseteq GA_j$. 

\begin{prop} \label{compact}
Assume that $\mathcal{B}$ is a computable nice basis 
corresponding to a compact $G$-space $({\bf X},\tau, G)$. 
Then there is a set ${\mathcal{O}}\subseteq \omega$ 
such that each $\varphi_e$ with $e\in {\mathcal{O}}$, 
codes a computable piece of the $\tau$-canonical 
partition which is a $G$-orbit, and all codes of 
computable closed $\tau$-canonical pieces 
which are $G$-orbits belong to ${\mathcal{O}}$. 
Moreover ${\mathcal{O}}$ belongs to $\Pi^0_3$ 
with respect to the complexity of $Con(z,U)$. 
\end{prop} 

{\em Proof.} 
Let ${\mathcal{O}}$ be the set of all $e$ 
satisfying $Con(e,{\bf X})$ such that 
for any $B\in {\mathcal{B}}$ one of the conditions 
$Con(e,GB)$ or $Con(e,{\bf X}\setminus GB)$ does not 
hold (i.e. $e$ codes a ${\bf t}$-canonical piece)
and for every $H\in {\mathcal{V}}^G$ there is 
a number $k$ such that for any $H$-invariant 
$C_1$,...,$C_{k+1}\in {\mathcal{B}}$ one of the 
conditions $Con(e, C_i \Delta C_j )$ does not hold.  
It is easy to see that ${\mathcal{O}}$ belongs 
to $\Pi^0_3$ with respect to complexity of 
$Con(z,U)$. 

As we have already mentioned above by Theorem 3.4 
of \cite{basia} the set ${\mathcal{O}}$ contains 
all codes of computable closed $\tau$-canonical pieces 
which are $G$-orbits. 
To see the proposition it remains to notice that 
if $e\in {\mathcal{O}}$, then the corresponding 
canonical piece $X_0$ has the property that 
for any $H\in {\mathcal{V}}^G$, any $H$-type 
$X_0$ is principal.  
Since there is only finitely many possibilities 
for intersections of $X_0$ with 
$H$-invariant members of ${\mathcal{B}}$ 
this claim is obvious. 
$\Box$
\bigskip 

{\bf Remark.} The case when the nice topology ${\bf t}$ 
is compact is not interesting. 
It does not differ from the case of the logic topology 
(i.e. logic $S_{\infty}$-space) of the first-order logic. 
In Proposition \ref{compact} the equality $\tau = {\bf t}$ 
corresponds to the latter case. 

\section{The automorphism group of a countably categorical structure} 

In this section we illustrate the material of Section 2 
in the case when the group $G$ is the automorphism group 
of an $\omega$-categorical structure with decidable theory. 
This slightly extends the corresponding material from \cite{LS} 
(where $G$ is $S_{\infty}$ and the topology is nice). 
We have found that the main construction 
of Section 2 of \cite{LS} is not presented in \cite{LS} in detail. 
Our Theorem 3.2 remedies this. 
Moreover it slightly generalizes the corresponding 
theorem of \cite{LS}.    

\subsection{Space} 

We fix a countable structure $M_0$ in a language $L_0$.
We assume that $M_0$ is $\omega$-categorical 
and the theory $Th(M_0)$ is decidable. 
Let $T$ be an extension of $Th(M_0)$ in a computable 
language $L$ with additional relational and functional 
symbols ${\bf r}_{1},...,{\bf r}_{t},...$ 
(possibly infintely many).
We assume that $T$ is axiomatizable
by first-order sentences of the following form:
$$
(\forall \bar{x})(\bigvee_{i} 
(\phi_{i}(\bar{x})\wedge\psi_{i}(\bar{x}))),
$$
where $\phi_{i}$ is a quantifier-free first-order formula 
in the language $L=L_0 \cup \{ {\bf r}_i \}_{i\in\omega}$, 
and $\psi_{i}$ is a first-order formula of the language $L_0$.
Consider the set $\mathbf{X}_{M_0}$ of all possible 
expansions of $M_0$ to models of $T$. \parskip0pt

For any tuple ${\bf \bar{r}}$ of ${\bf r}_i$-s and a tuple 
$\bar{a}\subset M_0$ we define as in \cite{ivanov} a {\em diagram} 
$\phi(\bar{a})$ of ${\bf \bar{r}}$ on $\bar{a}$.
To every functional symbol from ${\bf \bar{r}}$ we associate
a partial function from $\bar{a}$ to $\bar{a}$.
Choose a formula from every pair
$\{ {\bf r}_{i}(\bar{a}'),\neg {\bf r}_{i}(\bar{a}')\}$, where 
${\bf r}_{i}$ is a relational symbol from ${\bf \bar{r}}$ and 
$\bar{a}'$ is a tuple from $\bar{a}$ of the corresponding length.
Then $\phi(\bar{a})$ consists of the conjunction of the chosen 
formulas and the definition of the chosen functions
(so, in the functional case we look at $\phi(\bar{a})$ as
a tuple of partial maps). \parskip0pt

Consider the class ${\bf B}_{T}$ of all theories
$D(\bar{a})$, $\bar{a}\subset M_0$, such that each 
of them consists of $Th(M_0 ,\bar{a})$ and 
a diagram of some ${\bf \bar{r}}$ on $\bar{a}$ 
satisfied in some $(M_0 ,{\bf r}_i )_{i\in\omega}\models T$.
We order ${\bf B}_{T}$ by extension: $D(\bar{a})\le D'(\bar{b})$
if $\bar{a}$ consists of elements of $\bar{b}$ and 
$D'(\bar{b})$ implies $D(\bar{a})$ under $T$ 
(in particular, the partial functions defined in 
$D'$ extend the corresponding partial functions
defined in $D$).
Since $M_0$ is an atomic model, each element of 
${\bf B}_{T}$ is determined by a formula of the form
$\phi(\bar{a})\wedge\psi(\bar{a})$, where $\psi$ is 
a complete first-order formula for $M_0$ and $\phi$ 
is a diagram of some ${\bf \bar{r}}$ on $\bar{a}$.
The corresponding formula $\phi(\bar{x})\wedge\psi(\bar{x})$
will be called {\em basic}. \parskip0pt

On the set ${\bf X}_{M_0}$ 
of all $L$-expansions of the structure $M_0$ 
we consider the topology generated by basic open sets of the form
$Mod D(\bar{a})=\{ (M_0 ,{\bf r}'_i )_{i\in\omega}:$ 
$(M_0 ,{\bf r}'_i )_{i\in\omega} \models D(\bar{a})\}$,
$\bar{a}\subset M_0$.
It is easily seen that any $Mod D(\bar{a})$ is clopen.
We denote this basis by $\mathcal{A}$. 
The topology is metrizable: fix an enumeration
$(\bar{a}_{0},{\bf \bar{r}}_1 ), (\bar{a}_{1},{\bf \bar{r}}_2 ),...$ 
of $M^{< \omega}_0 \times (L\setminus L_{0})^{<\omega}$ and define 
\begin{quote}
$d((M_0 ,{\bf r}'_i )_{i\in\omega},(M_0 ,{\bf r}''_i )_{i\in\omega})$ 
$=\sum \{ 2^{-n}:$ 
there is a symbol ${\bf r}\in {\bf \bar{r}}_n$ such that 
its interpretations on $\bar{a}_{n}$ in the structures 
$(M_0 ,{\bf r}'_i )_{i\in\omega}$ and $(M_0 ,{\bf r}''_i )_{i\in\omega}$ 
are not the same (if ${\bf r}$ is a functional symbol then 
${\bf r'}_i (\bar{b}) \not= {\bf r''}_i (\bar{b})$ for some
$\bar{b}\subseteq \bar{a}_{n}$) $\}$.
\end{quote}

\noindent
It is easily seen that the metric $d$ defines the topology
determined by the sets of the form $Mod D(\bar{a})$. 
This topology will be denoted by $\mathbf{t}_{M_0}$.  
It is worth noting that by the assumptions on $T$ 
($T$ is axiomatizable by $\forall$-sentences  
with respect to symbols from ${\bf r}_i$) the space 
${\bf X}_{M_0}$ forms a closed subset of the space 
$\mathbf{X}_{L}$ of all $L$-structures on $\omega$. 
Thus ${\bf X}_{M_0}$ is a Polish space. 
\parskip0pt

Consider the action of the automorphism group  
$G:=Aut(M_0 )$ on the space $\mathbf{X}_{M_0}$. 
The basis $\mathcal{N}^G$ is defined to be all finite 
$Th(M_0 )$-elementary maps in $M_0$. 

\begin{lem} \label{niceB} 
The family of all sets $Mod(\phi (\bar{s}))$, where 
$\phi (\bar{s})$, $\bar{s}\in M_0$, is a first-order 
formula of the language $L$, 
forms a nice basis $\mathcal{B}$ of the $G$-space  
$(\mathbf{X}_{M_0} ,\mathbf{t}_{M_0})$.
\end{lem}  

{\em Proof.} 
This is verified in Theorem 1.10 of \cite{becker} for 
the $S_{\infty}$-space $\mathbf{X}_L$. 
Although the case of $\mathbf{X}_{M_0}$ is similar, 
some details are worth explaning. 
As in \cite{becker} we concentrate on condition (b)(iv) 
of the definition of a nice topology. 
We thus fix $B\in \mathcal{B}$ and $H\in \mathcal{N}^{G}$, and 
find pairwise distinct $r_0 ,...,r_{l-1}$, $s_0 ,...,s_{m-1}$, 
$t_0 ,...,t_{n-1}\in M_0$ and pairwise distinct 
$s'_0 ,...,s'_{m-1}$, $t'_0 ,...,t'_{n-1}\in M_0$ 
so that the following three conditions are satisfied: 

(1) the type of $s_0 ,...,s_{m-1}, t_0 ,...,t_{n-1}$ in $M_0$ 
coincides with the type of 

$s'_0 ,...,s'_{m-1}$, $t'_0 ,...,t'_{n-1}$;  

(2) $H= \{ g\in Aut(M_0 ): g(s'_0 )=s_0 ,...,g(s'_{m-1})=s_{m-1}$, 
$g(t'_0 )=t_0 ,...,g(t'_{n-1})$ 

$=t_{n-1}\}$; 

(3) $B= Mod (\phi (s_0 ,...,s_{m-1}, r_0 ,...,r_{l-1}) )$, 
where $\phi(\bar{u},\bar{z})$ is a first-order $L$-formula. \\ 
Let $\psi (u_0 ,...,u_{m-1},v_0 ,...,v_{n-1})$ be 
the following formula: 
$$
(\forall w_0 ,...,w_{l-1}) [(\mbox{ the type of } 
u_0 ,...,u_{m-1}, v_0 ,...,v_{n-1}, w_0 ,...,w_{l-1} \mbox{ in } M_0 
$$
$$
\mbox{ coincides with the type of }
s_0 ,...,s_{m-1}, t_0 ,...,t_{n-1}, r_0 ,...,r_{l-1}) \rightarrow 
$$
$$
\phi(u_0 ,...,u_{m-1}, w_0 ,...,w_{l-1})]. 
$$ 
Note that by $\omega$-categoricity of $M_0$, the first part 
of the implication above can be written by a first-order 
$L_0$-formula without parameters. 
To see that 
$$
B^{*H} = Mod(\psi , s'_0 ,...,s'_{m-1},t'_0 ,...,t'_{n-1})
$$ 
note that for any expansion $(M_0 ,{\bf r}'_i )$ satisfying 
$\psi (s'_0 ,...,s'_{m-1},t'_0 ,...,t'_{n-1})$, all automorphisms 
from $u$ take $(M_0 ,{\bf r}'_i )$ to $B$. 
On the other hand if the expansion $(M_0 ,{\bf r}'_i )$ 
does not satisfy $\psi (s'_0 ,...,s'_{m-1},t'_0 ,...,t'_{n-1})$, 
then there is a tuple $r'_0 ,...,r'_{l-1}$ such that 
the basic open set of all automorphisms of $M_0$ 
defined by the map 
$$
s'_0 ,...,s'_{m-1}, t'_0 ,...,t'_{n-1}, r'_0 ,...,r'_{l-1}\rightarrow 
s_0 ,...,s_{m-1}, t_0 ,...,t_{n-1}, r_0 ,...,r_{l-1}
$$ 
is non-empty and does not contain an element 
taking $(M_0 ,{\bf r}'_i )$ to $B$. 
$\Box$
\bigskip 

To check that the $G$-space ${\bf X}_{M_0}$ satisfies 
the computability conditions above, note that 
$M_0$ has a presentation on $\omega$ so that all relations 
first-order definable in $M_0$, are decidable. 
This follows from $\omega$-categoricity and decidability 
of $Th(M_0 )$ together with the standard fact that 
a decidable theory has a strongly constuctivizable model. 
We fix such a presentation. 
Then we can define a computable presentation of the 
following sorts and relations: 
the elements of $\mathcal{V}^{G}$ can be interpreted by 
finite subsets of $M_0$ and elements of $\mathcal{N}^{G}$ 
are interpreted by elementary functions between finite subsets 
of $M_0$. 
Since the elementary diagram of $M_0$ is decidable, 
the set of elementary functions between finite subsets of $M_0$ is computable. 

We can also consider elements of $\mathcal{V}^{G}$ as 
finite identity functions. 
The relation of inclusion $\subset$ on $\mathcal{N}^G$ is 
defined by $g_1 \subseteq g_2 \Leftrightarrow$ "$g_2$ is 
a restriction of $g_1$". 
When we consider elements of $\mathcal{V}^{G}$ as
finite identity functions, this inclusion corresponds to 
the standard one on $\mathcal{V}^G$.
\parskip0pt 

Since we interpret elements of $\mathcal{B}$ by $L$-formulas 
with parameters from $M_0$ and without free variables,
it is obvious that $\mathcal{B}$ can be coded in $\omega$ 
so that the operations of the Boolean algebra $\mathcal{B}$  
are defined by decidable predicates. 
For example the operations $\neg$, $\wedge$ and $\vee$ play 
the role of $'$, $\cap$ and $\cup$. 
The operation of taking $*$-transform is coded according 
the construction of the proof of Lemma \ref{niceB}.
Since the basis $\mathcal{A}$ is interpreted by quantifier-free 
formulas, it is a decidable subset of $\mathcal{B}$. 
Then $\cap$ and $\cup$ define the ordering of $\mathcal{A}$. 
The remaining basic relations are defined as follows. \parskip0pt

$Inv(V,U)\Leftrightarrow$ 
"the parameners of $U$ are uniquely defined in $M_0$ over the set $V$"
(i.e. if $U$ is a basic subset defined by an $L$-formula 
$\phi$ with parameters $\bar{a}$ and $V$ is the $G$-stabiliser 
of a tuple $\bar{c}$, then there is an $L_0$-formula 
$\psi (\bar{x},\bar{c})$ over $\bar{c}$  such that 
$M_0 \models\forall\bar{x}(\psi(\bar{x},\bar{c})\rightarrow\bar{x}=\bar{a})$); 

$Orb_{l,n} (N,V_1 ,...,V_l ,V_{l+1},...,V_{2l} ,U_{1},...U_n ,U_{n+1},...,U_{2n})\Leftrightarrow$ 
$N\in \mathcal{N}^{G}\wedge $ 

$\bigwedge^{2m}_{i=1}(V_i \in \mathcal{V}^{G})\wedge$ 
$\bigwedge^{2m}_{i=1}(U_i \in \mathcal{A})\wedge$
"there is an $M_0$-elementary bijection $g$ 

between the set of all elements arising as stabilized 
points of $V_1 ,...,V_m$ 

and/or as parameters of the formulas $U_1 ,...,U_n$ and 
the corresponding set 

arising in $V_{m+1},...,V_{2m}$ and the formulas 
$U_{n+1},...,U_{2n}$ such that $g$ extends 

the map defining $N$ and maps each $V_i$ (the code 
of each $U_i$) to $V_{i+m}$, 

$i\le m$ (to the code of $U_{n+i}$, $i\le n$)". \\
By $\omega$-categoricity and decidability of the chosen presentation 
of $M_0$, these relations are also decidable. 
 
Let $\phi(\bar{s})$ be a quantifier-free formula defining an element 
$A\in \mathcal{A}$.  
To compute $diam(A)$ consider the definition of the metric $d$ above. 
Using decidability of the elementary diagram of $M_0$ find the greatest $n$ 
such that for all $i\le n$ the interpretation of ${\bf \bar{r}}_i$ on 
$\bar{a}_i$ is uniquely determined by $\phi (\bar{s})$. 
Then $2^{-n-1}\le diam(A)<2^{-n}$. 
The case of basic clopen sets of $\mathcal{N}^G$ is similar.

\subsection{Examples} 

In the case of ${\bf X}_{M_0}$ we can use the argument 
of Section 2 of \cite{LS} to show that the class 
$\Pi^0_3$ contains the set of all numbers of 
${\bf t}_{M_0}$-canonical pieces, which are $G$-orbits. 
To see this note that each canonical piece is defined by 
sentences of the form $\exists \bar{x} D(\bar{x})$ and 
$\neg \exists \bar{x} D(\bar{x})$, where $D(\bar{x})$ is a basic formula. 
If the corresponding theory of such sentences together with 
$Th(M_0 )$ axiomatizes an $\omega$-categorical $L$-theory, 
then the canonical piece is a $G$-orbit. 
When the corresponding theory is not $\omega$-categorical 
then by $\omega$-categoricity of $Th(M_0 )$ we can find 
two $L$-expansions of $M_0$ of our canonical piece which are not 
isomorphic, i.e. are not in the same $G$-orbit. 

We now see that to state that a canonical piece of 
${\bf X}_{M_0}$ is a $G$-orbit it is enough to express that 
the corresponding $L$-theory (together with $Th(M_0 )$) 
satisfies the conditions of the Ryll-Nardzewski theorem  
(i.e. we have finitely many $n$-types for all $n$). 
It is shown in \cite{LS} that this can be written as 
a $\Pi^0_3$-condition. 
The following theorem roughly claims that the set of 
canonical pieces which are $G$-orbits, is $\Pi^0_3$-complete. 

\begin{thm} 
Let $N$ be an $\omega$-categorical infinite 
structure with decidable theory. 
Then there is a decidable $\omega$-categorical 
(say $L_0$)-structure $M_0$ such that $N$ is 
interpreted in $M_0$ and for some infinite 
language $L\supset L_0$ there is an $L$-theory 
$T$ extending $Th(M_0 )$ and satisfying the assumptions  
of Section 3.1 (in particular $\forall$-axiomatizability 
with respect to $L\setminus L_0$) 
such that the $Aut(M_0 )$-space $\mathbf{X}_{M_0}$ of 
the $L$-expansions has the canonical partition with 
the property that the set of 
all natural numbers $e$ satisfying the relation 
 
"$\varphi_e$ codes a piece of the canonical partition which is 
an $Aut(M_0 )$-orbit"  \\  
is $\Pi^{0}_3$-complete. 
\end{thm} 

{\em Proof.} 
The proof is based on two constructions: 
\begin{quote} 
$*$ the idea of Section 2 of \cite{LS} of the proof for 
the case when $M_0$ is a pure set; \\
$*$ the construction of $\omega$-categorical expansions from \cite{P}. 
\end{quote} 
We start with the presentation of the latter one. 
Let $L_E$ consist of $2n$-ary relational symbols $E_n$, 
$n\in\omega\setminus \{ 0\}$, and $T_E$ be the 
$\forall\exists$-theory of the universal homogeneous 
structure of the universal theory saying that each $E_n$ is 
an equivalence relation on the set of $n$-tuples 
such that all $n$-tuples with at least one repeated 
coordinate lie in one isolated $E_n$-class.  \parskip0pt 

Let $T'$ be a many-sorted $\omega$-categorical 
theory in a relational language $L'$ 
with countably many sorts $S_n$, $n\in\omega$, 
such that elements of $S_0$ may appear only in $=$. 
Let $M$ be a countable model of $T_E$ and $M_{\bar{S}}$ 
be the expansion of $M$ to the language 
$L_E \cup \{ S_1, ...,S_n ,...\} \cup \{ \pi_1 ,..., \pi_n ,...\}$, 
where each $S_n$ is interpreted by the non-diagonal elements 
of $M^n /E_n$ and $\pi_n$ by the corresponding projection. 
By $(M_{\bar{S}})'$ we denote a $T'$-expansion of $M_{\bar{S}}$ 
to the language $L'$, where $S_0$ is identified with the 
basic sort of $M$. 
Theorem 4.2.6 of \cite{P} states that all such expansions 
have the same theory and this theory is $\omega$-categorical.  
\parskip0pt 

We now build an expansion $M^*$ of $M$ (in the 1-sorted language). 
For each relational symbol $R_i \in L'$ of the sort 
$S_{n_1}\times S_{n_2}\times ...\times S_{n_k}$ we 
add a new relational symbol $R^*_i$ on 
$M^{n_1 \cdot ...\cdot n_k}$ interpreted in the following way: 
$$
M^* \models R^*_i (\bar{a}_1 ,...,\bar{a}_k) \Leftrightarrow 
(M_{\bar{S}})'\models R_i (\pi_{n_1}(\bar{a}_1 ),...,\pi_{n_k}(\bar{a}_k )). 
$$  
It is clear that $M^*$ and $(M_{\bar{S}})'$ are bi-interpretable. 
Thus $Th(M^* )$ is $\omega$-categorical. 

We now prove the main statement of the theorem. 
Let $N$ be an $\omega$-categorical structure. 
Let $L_0$ be $L_E$ together with the language of $N$ 
(where the basic sort is denoted by $S_0$ as above). 
To define $L_1$, for every natural $n\ge 2$ we extend 
$L_0 \cup \{ S_1 ,...,S_n ,...\} \cup \{ \pi_1 ,...,\pi_n ,...\}$ 
by an $\omega$-sequence of unary relations $P_{n,i}$, 
$i\in\omega$, defined on $S_n$. 
We also put all relations of $N$ onto the sort $S_1$.   
Let $T_1$ be the $L_1$-theory axiomatized by $T_E$ 
together with the natural axioms for 
all $\pi_n$, with the theory $Th(N)$ on $S_1$ and with 
the axioms saying that all $N$-relations on $S_0$ are 
just $*$-versions of $N$-relations on $S_1$. 
By $T$ we denote the theory of all $M^*$ with $(M_{\bar{S}})'\models T_1$. 
Let $L$ be the corresponding language. 
Let $M_0$ be the $L_0$-reduct of a countable $M^*\models T$. 
It is clear that $T$ is axiomatized by $Th (M_0 )$ (containing  
the $\forall\exists$-axioms of $T_E$) and $\forall$-axioms of 
$E_n$-invariantness of $P^*_{n,i}$, $n\ge 2$, $i\in\omega$.  
Thus $M_0$ and $T$ satisfy the basic assumptions 
of the previous subsection. 
In particular $Th(M_0 )$ is $\omega$-categorical and 
decidable by Theorem 4.2.6 of \cite{P} (cited above) and 
by $\omega$-categoricity and decidability of $Th(N)$ 
(the latter implies that $Th(M_0 )$ is computably axiomatizable). 
\parskip0pt 

For every sequence of finite sets of natural numbers 
$\theta =(D_2 ,D_3, ...,D_n ,...)$ we define the many-sorted 
$L_1$-theory $T_{\theta}\supset T_1$ saying that for each $n$, all 
$P_{n,j}$ with $j\not\in D_n$, are empty, and the family $P_{n,j}$, 
$j\in D_{n}$ freely generates a Boolean algebra of infinite subsets 
of $S_n$ (denote the $n$-th part of $T_{\theta}$ by $T_{n,D_n}$). 
Again by Theorem 4.2.6 of \cite{P} each $T_{\theta}$ is 
$\omega$-categorical. 
Moreover it is obtained from $T_1$ by adding some axioms 
which are just $\forall$- or $\exists$-sentences concerning $P_{n,j}$. 
\parskip0pt 

Let $M$ be a countable $L_0$-model of $Th (M_0 )$. 
By $M_{\theta}$ we denote an expansion of $M$ to $T_{\theta}$. 
As we already know, by Theorem 4.2.6 of \cite{P}, all these 
expansions are $\omega$-categorical and isomorphic.   
Since they are axiomatized by $T_E$, $Th(N)$ (on $S_1$) 
and all $T_{n,D_n}$, $n\in\omega$, we see that for any 
two sequences $\theta' =(D'_2 ,D'_3 ,...,D'_n ,...)$ 
and $\theta'' =(D''_2 ,D''_3 ,...,D''_n ,...)$ with 
$D''_n \subset D_n \cap D'_n$, $n\in\omega$, the reducts 
of $M_{\theta}$ and $M_{\theta'}$ to 
$L_0 \cup \bigcup \{ P^*_{n,i}: i\in D''_n , n\in\omega\}$
are isomorphic. 

For every natural $e$ let us fix a computable enumeration 
$\rho_e$ (as a function defined on $\omega$) of the set of 
all pairs $\langle n,x\rangle$ with $x\in W_{\varphi_{e}(n)}$. 
For every natural $l$ we define a sequence 
$\theta_l =(D_2 ,D_3 ,...)$ of finite sets such that 
$$
k\in D_n \Leftrightarrow (k\le l)\wedge 
(\exists x)(\rho_e (k)=\langle n-2 ,x\rangle \wedge 
(\forall k'<k)(\rho_e (k')\not=\langle n-2 ,x\rangle )).
$$
Let $T_e$ be the $L_1$-theory such that for every natural 
$l$ the reduct of $T_e$ to 
$$
L_0 \cup \{ S_1 ,...,S_n ,...\} \cup \{ \pi_1 ,..., \pi_n ,...\} 
\cup\{ P_{n,i}, i\le l \mbox{ and } 2\le n\}
$$ 
coincides with the corresponding reduct of $T_{\theta_l }$. 
It is obvious that $T_e$ is axiomatizable by 
a computable set of axioms (uniformly in $e$). 
Since for each $l$ the reduct of $T_e$ as above is 
$\omega$-categorical, the theory $T_e$ is complete. 
Thus $T_e$ is decidable uniformly in $e$.  
By Ryll-Nardzewski's theorem the theory $T_e$ is 
$\omega$-categorical if and only if all $W_{\varphi_{e}(k)}$ 
are finite (i.e. the set of 1-types of each $S_k$ is finite). 
If we consider models of $T_e$ in the 1-sorted $*$-form 
defined as above, then these properties remain true. 

Let $M_0$ be as above.
As we have already mentioned $M_0$ is $\omega$-categorical, 
the theory $Th(M_0)$ is decidable and the theory $T$ is 
an $L$-extension of $Th(M_0)$ which is axiomatizable
by first-order sentences of the following form:
$$
(\forall \bar{x})(\bigvee_{i} 
(\phi_{i}(\bar{x})\wedge\psi_{i}(\bar{x}))),
$$
where $\phi_{i}$ is a quantifier-free first-order formula 
in the language $L$ and $\psi_{i}$ is a first-order formula 
of the language $L_0$.
Consider the space $\mathbf{X}_{M_0}$ of all possible 
expansions of $M_0$ to models of $T$. 
The group $G=Aut(M_0 )$ makes it a Polish G-space. 
\parskip0pt

Since the $*$-form of each $T_e$ is a decidable complete theory 
axiomatized by $Th(M_0 )$ and universal/existentional sentences 
concerning all $P^{*}_{n,i}$, all the structures of $\mathbf{X}_{M_0}$ 
corresponding to $T_e$ form a computable piece of the canonical 
partition on $\mathbf{X}_{M_0}$.  
Since any algorithm computing $\varphi_e$ effectively provides 
an algorithm deciding the $*$-version of $T_e$ with respect to 
existential/universal $P^{*}_{n,i}$-sentences, we easily see that 
the $\Pi^{0}_{3}$-set $\{ e: \forall n (W_{\varphi_e (n)}$ is finite)$\}$ 
is reducible to $\{ e: T_e$ is $\omega$-categorical$\}$.  
Since the former one is $\Pi^{0}_{3}$-complete (see \cite{LS} 
and \cite{soare}, p.68) we have the theorem. $\Box$

\bigskip 

{\bf Remark.} 
Analysing examples of \cite{ivanov} and \cite{P} one can 
prove that the statement of the theorem holds for the class 
$\Pi^0_2$ and the relation 
\begin{quote} 
"$\phi_e$ codes a piece of the canonical partition which 
is an $Aut(M_0 )$-orbit of a G-compact structure". 
\end{quote}
The definition of G-compacness can be also found 
in \cite{ivanov} and \cite{P}. 
Since this notion is not so natural outside model theory, 
we do not develop this further.  

\section{Degree spectrum of canonical pieces} 

\subsection{The space ${\bf X}_{M_0}$}

In this section we preserve the assumptions of Section 2. 
Let $G$ be a closed subgroup of $S_{\infty}$ and 
$({\bf X},\tau )$ be a Polish $G$-space. 
Let $\mathcal{A}$ be a countable basis of 
$({\bf X},\tau )$ closed with respect to $\cap$. 
Each $A\in \mathcal{A}$ is $H$-invariant with
respect to some basic subgroup $H\in \mathcal{N}^G$. 
The subfamily of $\mathcal{A}$ consisting of clopen 
sets generates the same topology. 
The bases $\mathcal{N}^G$ and $\mathcal{A}$ 
are computably 1-1-enumerated so that the relations 
$\subseteq$, $\cap$, $Clopen$,  $Inv(V,U)$ and 
$$
Orb_{m,n}(N,V_1 ,...,V_m ,V_{m+1},...,V_{2m},U_1 ,...,U_n ,U_{n+1},...,U_{2n})
$$
are presented by decidable relations on $\omega$. 
There is an algorithm deciding 
the problem if for a basic clopen set $U$ 
(of $\mathcal{N}^G$ or $\mathcal{A}$) and a natural 
number $i$ the diametr of $U$ is less than $2^{-i}$.  

\begin{definicja} \label{Degree} 
We say that an element $x\in {\bf X}$ {\em represents degree unsolvability} 
${\bf d}$ if the relation 
$$
Sat_x (U)\Leftrightarrow (U\in \mathcal{A})\wedge (x\in U)
$$ 
(i.e. the set $Sat_x ({\mathcal{A}})$) is of degree ${\bf d}$. 
\end{definicja} 

In the case of the logic action, when $x$ is a structure 
on $\omega$, this notion is obviously equivalent to 
the notion of a structure of degree ${\bf d}$.  
As before it is straightforward that for an $x$ of degree ${\bf d}$ 
there is a ${\bf d}$-comutable function 
$\kappa :\omega \rightarrow \mathcal{A}$ 
such that for all $n$, $x\in \kappa (n)$ 
and $\kappa (n)$ is clopen with $diam(\kappa (n))<2^{-n}$. 
It is also worth noting that when $\mathcal{A}$ consists of 
clopen sets the existence of such ${\bf d}$-computable $\kappa$ 
already implies that the set 
$Sat_x ({\mathcal{A}})$ is of degree ${\bf d}$. 

We say that the {\em orbit $Gx$ is of degree} ${\bf d}$ if 
${\bf d}$ is the least degree of the members of $Gx$.  
In the case when such a degree does not exist we say that 
$Gx$ {\em has no degree}. 

Following \cite{richter} we now introduce combination methods 
for $\mathcal{A}$. 
We say that a computable subfamily $A_1 ,...,A_n ,...$ of $\mathcal{A}$ 
is {\em effectively free} if every its finite subfamily freely generates 
a Boolean algebra of sets.  
The following theorem is a counterpart of Theorem 2.1 of \cite{richter}. 

\begin{thm} 
Let $A_1 ,...,A_n ,...$ be an effectively free subfamily of $\mathcal{A}$. 
Assume that for each $S\subseteq \omega$ there exists 
an element $x_{S}\in {\bf X}$ such that \\ 
(i) $Sat_{x_{S}}({\mathcal{A}})$ is computable with respect 
to $S$ and \\ 
(ii) $\forall i\in \omega ( Sat_{x_{S}} (A_i ) \Leftrightarrow i\in S)$. 

Then for every degree ${\bf d}$ there is an element $x\in {\bf X}$ 
such that the orbit $Gx$ is of degree ${\bf d}$. 
\end{thm} 

{\em Proof.} A straightforward adaptation of the proof 
of Theorem 2.1 from \cite{richter}. $\Box$

\bigskip 

We now consider the case when $Gx$ has no degree. 

\begin{thm} 
Let $A_1 ,...,A_n ,...$ be an effectively free subfamily of $\mathcal{A}$. 
Assume that for each $S\subseteq \omega$ there exists 
an element $x_{S}\in {\bf X}$ such that \\ 
(i) $Sat_{x_S} ({\mathcal{A}})$ is enumeration reducible to $S$ 
\footnote{there is an effective procedure whose outputs enumerate 
$Sat_{x_S}({\mathcal{A}})$ when any enumeration of $S$ is supplied 
for the inputs} 
and \\ 
(ii) $\forall i\in\omega ( Sat_{x_{S}} (A_i ) \Leftrightarrow i\in S )$. 

Then there is a set $S$ such that the orbit $Gx_{S}$ has no degree. 
\end{thm} 

{\em Proof.} A straightforward adaptation of the proof 
of Theorem 2.3 from \cite{richter}. 
We just remind the reader that it is based on the fact that 
there exists a set $S\subset \omega$ such that the mass problem 
$\{ f: range(f)=S\}$ has no Turing-least element. 
Having such an $S$ it is straightforward to show that 
the mass problem ${\bf E}_{S}= \{ f: range(f)=S\}$ is 
Medvedev-equivalent to the problem ${\bf Ch}_{Gx_S}$ 
of all characteristic functions of all sets 
$Sat_{x}({\mathcal{A}})$, $x\in Gx_{S}$. 
This means that there are partial computable operators 
$\Phi$ and $\Psi$ such that $\Phi$ maps ${\bf E}_S$ 
to ${\bf Ch}_{Gx_S}$ and $\Psi$ maps ${\bf Ch}_{Gx_S}$ 
to ${\bf E}_S$. 
Since for total functions Turing-reducibility coincides 
with the enumeration reducibility (see Chapter 9 of \cite{rogers}) 
the existence of the least Turing degree of $Gx_S$ 
(i.e. of ${\bf Ch}_{Gx_S}$) implies the same property 
for ${\bf E}_{S}$, a contradiction. $\Box$

\bigskip 

We can now present the main results of this section. 

\begin{thm} 
Let $N$ be an $\omega$-categorical model complete 
infinite structure with decidable theory. 
Then the decidable $\omega$-categorical 
$L_0$-structure $M_0$ (such that $N$ is interpreted in $M_0$), 
the infinite language $L\supset L_0$ and the $L$-theory $T$ 
(extending $Th(M_0 )$) constructed in Theorem 3.2 have the 
property that the canonical partition of the $Aut(M_0 )$-space 
$\mathbf{X}_{M_0}$ of the $L$-expansions has \\ 
(i) canonical pieces which are G-orbits of any possible degree ${\bf d}$; \\ 
(ii) canonical pieces which are G-orbits having no degree. 
\end{thm} 

{\em Proof.} 
We now apply the construction of the proof of Theorem 3.2. 
Let $L_E$, $L_0$ and $L$ be as in that proof.    
We also repeat the definition of $T_E$, $T_1$ and $T$ 
(the theory of all $M^*$ with $M\models T_1$). 
As above $M_0$ is the $L_0$-reduct of a countable $M^*\models T$. 
Fix any computable enumeration of $T$. 
\parskip0pt 

In the proof of Theorem 3.2 for every sequence of finite sets 
of natural numbers $\theta =(D_2 ,D_3, ...,D_n ,...)$ we have 
defined the many-sorted $\omega$-categorical theory 
$T_{\theta}\supset T_1$ saying that for each $n$, all $P_{n,j}$ 
with $j\not\in D_n$, are empty, and the family $P_{n,j}$, 
$j\in D_{n}$, freely generates a Boolean algebra of infinite 
subsets of $S_n$ (where the $n$-th part of $T_{\theta}$ is 
denoted by $T_{n,D_n}$). 
\parskip0pt 

For a subset $S\subseteq \omega$ by $M_S$ we denote the expansion 
$M_{\theta}\models T_{\theta}$, where $\theta =(D_2 ,...,D_n ,...)$ 
with $D_{i+2} =\{ 1\}$ for $i\in S$, and 
$D_{i+2} = \emptyset$ for $i\not\in S$. 
It is clear that each $(M_S )^*$ is $\omega$-categorical.  
Since $Th(N)$ is model complete, the theory $Th((M_S )^* )$ 
is $\forall\exists$-axiomatizable and thus model complete too. 
Since its axioms are computable in $S$, it is decidable in $S$. 
In particular $(M_S)^*$ has a presentation such that its 
elementary diagram is computable in $S$. 

Any enumeration of $S$ provides an enumeration 
of an infinite substructure of $(M_S )^*$ as follows. 
Assume that at step $n-1$ we have already 
enumerated a subset $Q\subset (M_S )^*$. 
Take the $n$-th initial segment of $S$ and find 
the maximal element $m$ in it. 
Consider all quantifier free formulas of the form 
$\phi (q_1 ,...,q_l ,x_1 ,...,x_k )$ with $q_i \in Q$, 
$0<k\le m_2$ and $0\le l$, which appear in the $n$-th 
initial segment of the enumeration of axioms 
of $T$ of the form 
$\forall z_1 ,...,z_l \exists x_1 ,....,x_k \phi (\bar{z},\bar{x})$ 
and in additional axioms of 
$Th((M_S )^* )$ of the form 
$\exists x_1 ,...,x_k \phi(\bar{x})$.   
Choosing some realizations of each formula of this form  
we extend $Q$ by these realizations. 
By categoricity and model completeness this procedure 
gives a structure isomorphic to $(M_S )^*$. 

Now consider the space $\mathbf{X}_{M_0}$ of all 
possible expansions of $M_0$ to models of $T$. 
The group $G=Aut(M_0 )$ makes it a Polish G-space. 
Moreover the $G$-orbit of $(M_S )^*$ as above is a piece of 
the canonical partition. 
Let $A_i$ be the basic set of all $T$-structures on $\omega$ 
which satisfy the elementary diagram $D^{M_0}(1,...,i+2)$ 
of the tuple $(1,...,i+2)$ in $M_0$ together 
with $P_{i+2,1}(1,....,i+2)$. 
By the definition of $T$ and Theorem 4.2.6  
of \cite{P} for every sequence $\varepsilon_i \in \{ 0,1 \}$, 
$i\le l$, the  formula of the form 
$\bigwedge_{i\le l} (D^{M_0}(1,...,i+2) \wedge P^{\varepsilon_i}(1,...,i+2))$ 
is realized by a $T$-structure on $\omega$.  
We conclude that the sequence $A_1 ,...,A_n ,...$ is 
an effectively free subfamily of the standard basis of 
$\mathbf{X}_{M_0}$ (defined by all diagrams as in Section 3.1). 
Now for every subset $S$ of $\omega$ the structure $(M_S )^*$
as above satisfies the conditions of Theorems 4.2 and 4.3. 
Note that condition (i) of each of these theorems easily follows 
from the properties of $(M_S )^*$ mentioned above. 
For example the enumeration constructed in the previous paragraph 
easily gives an enumeration of $Sat_{(M_S )^*}({\mathcal{A}})$. 
This proves our theorem. $\Box$

\subsection{Countably categorical groups} 

It is worth noting that the construction of 
the previous subsection also gives examples of 
structures such that their isomorphism types 
have (have no) degree. 
Since these structures are $\omega$-categorical 
it seems to the authors that the examples are really new. 
In particular they provide theories having (having no) degrees. 

Sometimes it is interesting to verify if examples of this 
kind can be found in natural algebraic classes: see \cite{DDHS} 
and \cite{HKSS}.  
In this section we consider $\omega$-categorical 2-step 
nilpotent groups with quantifier elimination. 
Using \cite{CSW} we give a construction of new examples.

We start with a description of a QE-group of nilpotency 
class 2 given in \cite{CSW}. 
Since the group is built as the Fra\"{i}ss\'{e} limit of a class of 
finite groups, we give some standard preliminaries 
(see for example \cite{evans}). \parskip0pt 
   
Let $\mathcal{K}$ be a non-empty class of finite structures of 
some finite language $L$. 
We assume that $\mathcal{K}$ is closed under isomorphism and under 
taking substructures (satisfies HP, the {\em hereditary property}), 
has the {\em joint embedding property} (JEP) and 
the {\em amalgamation property} (AP). 
The latter is defined as follows: for every pair of embeddings  
$e: A \rightarrow B$ and $f: A \rightarrow C$ with 
$A,B,C\in \mathcal{K}$ there are embeddings $g: B \rightarrow D$ 
and $h: C \rightarrow D$ with $D\in \mathcal{K}$ such that 
$g\cdot e =h\cdot f$. 
Fra\"{i}ss\'{e} has proved that under these assumptions there 
is a countable locally finite 
\footnote{i.e. every finitely generated substructure is finite} 
$L$-structure $M$ (which is unique up to isomorphism) such that: 
\begin{quote}
(a) $\mathcal{K}$ is the {\em age} of $M$, i.e. the class of all 
finite substructures which can be embedded into $M$ and \\
(b) $M$ is {\em finitely homogeneous} (ultrahomogeneous), 
i.e. every isomorphism between finite substructures of $M$ 
extends to an automorphism of $M$.   
\end{quote}
The structure $M$ is called the {\em Fra\"{i}ss\'{e} limit} of $\mathcal{K}$. 
It admits elemination of quantifiers.  
\parskip0pt 

To define a $2$-step nilpotent, $\omega$-categorical 
homogeneous groups we assume that $\mathcal{K}$ is 
the class of all finite groups of exponent four in which 
all involutions are central.  
By \cite{CSW} $\mathcal{K}$ satisfies the HP, the JEP and the AP. 
Let $\mathcal{G}$ be the Fra\"{i}ss\'{e} limit of this class. 
Then $\mathcal{G}$ is nilpotent of class two. \parskip0pt 

We need the notions of free amalgamation and 
a-indecomposability in $\mathcal{K}$. 
Following \cite{CSW} we define them through the associated 
category of {\em quadratic structures}. 
A quadratic structure is a structure $(U,V;Q)$ where $U$ and 
$V$ are vector spaces over the field ${\bf F}_2$ and 
$Q$ is a nondegenerate quadratic map from $U$ to $V$, i.e. 
$Q(x) \neq 0$ for all $x \neq 0$ and the function 
$\gamma(x,y)=Q(x)+Q(y)+Q(x+y)$ is an alternating 
bilinear map. 
By $\mathcal{Q}$ we denote the category of all 
quadratic structures with morphisms 
$(f,g):(U_1,V_1;Q_1) \rightarrow (U_2,V_2;Q_2)$
given by linear maps $f:U_1 \rightarrow U_2$, 
$g:V_1 \rightarrow V_2$ respecting the quadratic 
map: $g Q_1 =Q_2 f$. \parskip0pt 

For $G \in \mathcal{K}$ define $V(G):= \Omega (G)$, 
the subgroup of all involutions of $G$, and $U(G):= G/V(G)$. 
Let $Q_G: U(G)\rightarrow V(G)$ be the map induced by squaring in $G$. 
Then $QS(G)= (U(G),V(G);Q_G)$ is a quadratic structure 
and the associated map $\gamma (x,y)$ is the one induced by 
the commutation from $G/V(G) \times G/V(G) $ to $V(G)$. 
It is shown in Lemma 1 of \cite{CSW} that this gives 
a 1-1-correspondence between $\mathcal{K}$ and $\mathcal{Q}$ 
up to the equivalence of central extensions 
$1\rightarrow V(G)\rightarrow G\rightarrow U(G)\rightarrow 1$ 
with $G\in \mathcal{K}$. \parskip0pt    

We now consider the amalgamation process in $\mathcal{K}$. 
To any amalgamation diagram in $\mathcal{K}$, 
$G_0 \rightarrow G_1 ,G_2$ we associate the diagram 
$QS(G_0 )\rightarrow QS(G_1 ),QS(G_2 )$ of the corresponding 
quadratic structures and (straightforward) morphisms.   
Let $QS(G_i )=(U_i ,V_i ;Q_i )$, $i\le 2$. 
Let $U^*, V^*$ be the amalgamated direct sums 
$U_1 \bigoplus_{U_{0}} U_2$, $V_1 \bigoplus_{V_{0}} V_2$ 
in the category of vector spaces. 
We define the free amalgam of $QS(G_1 )$ and $QS(G_2 )$ over 
$QS(G_0 )$ as a quadratic structure $(U ,V ;Q)$ with $U=U^*$ and 
$V=V^* \bigoplus (U_1 /U_0 )\otimes (U_2 /U_0 )$ (see \cite{CSW}). 
The corresponding quadratic map $Q:U \rightarrow V$ 
is defined by first choosing splittings of $U_1$, $U_2$ 
as $U_0 \bigoplus U'_1$ and $U_0 \bigoplus U'_2$, 
respectively, identifying $U'_1$, $U'_2$ 
with $U_1 /U_0$, $U_2 /U_0$ and defining 
$$
Q(u_0 +u'_1 +u'_2 )= Q_0 (u_0) +Q_1 (u'_1) +Q_2 (u'_2) + 
\gamma_1 (u_0, u'_1 ) + \gamma_2 (u_0,u'_2 ) + (u'_1 \otimes u'_2 ).
$$
Note that $Q|_{U_i} =Q_i$
%$Q \upharpoonright U_i =Q_i $ 
and the corresponding 
$\gamma (u'_1, u'_2)$ is $u'_1 \otimes u'_2$. 
Since $u'_1 \otimes u'_2 =0$ only when one of the factors 
is zero, the nondegeneracy is immediate.
It is shown in \cite{CSW} that $(V,U;Q)$ is a pushout of 
the natural maps $QS(G_1 )$, $QS(G_2 )\rightarrow (V,U;Q)$ 
agreeing on $QS(G_0 )$. 
We call the quadratic structure $(V,U;Q)$ the {\em free amalgam} 
of $QS(G_1 )$, $QS(G_2 )$ over $QS(G_0 )$. 
Let $G$ be the group associated with $(V,U;Q)$ in $\mathcal{K}$. 
By Lemma 3 of \cite{CSW} there are embeddings 
$G_1, G_2 \rightarrow G$ with respect to which $G$ becomes 
an amalgam of $G_1$, $G_2$ over $G_0$ in $\mathcal{K}$. 
We call $G$ {\em the free amalgam} of $G_0 \rightarrow G_1 ,G_2$. 
\parskip0pt 

We call a group $H\in \mathcal{K}$ {\em a-indecomposable} 
if whenever $H$ embeds into the free amalgam of two structures 
over a third, the image of the embedding is contained in one 
of the two factors. 
It is proved in Section 3 of \cite{CSW} that there is 
a sequence of a-indecomposable groups 
$\{ G_d :d\in\omega\} \subseteq \mathcal{K}$ such that for 
any pair $d\not= d'$ the group $G_d$ is not embeddable into $G_{d'}$.  
The construction is as follows. 
For any prime $p$ let $\hat{F}_p =(GF(2^{2p}), GF(2^p ); N)$ be 
the quadratic structure consisting of the finite fields of 
orders $2^{2p}$ and $2^p$ respectively and the corresponding 
norm $N:GF(2^{2p})\rightarrow GF(2^p )$. 
By Lemmas 9 and 12 of \cite{CSW} the sequence of the 2-step 
nilpotent groups $G_n$, $n\in \omega$, corresponding to  
the quadratic structures $\hat{F}_{p_n}$, $n\in \omega$, 
gives an appropriate antichain.  

It is worth noting that the construction is effective 
in the following sense. 
Since ${\mathcal{K}}$ consists of finite 
structures, we find an effective enumeration of ${\mathcal{K}}$ 
by natural numbers. 
Then the set of all groups $G_n$ forms a computable subset 
of the class $\mathcal{K}$.  

\begin{thm} 
(1) For every degree ${\bf d}$ there is an $\omega$-categorical 
2-step nilpotent QE-group $G$ of exponent four 
such that the isomorphism class of $G$ is of degree ${\bf d}$. 

(2) There is an $\omega$-categorical 2-step nilpotent 
QE-group $G$ of exponent four such that the isomorphism 
class of $G$ has no degree. 
\end{thm} 

{\em Proof.} (1) We apply Theorem 2.1 from \cite{richter} 
to the effective antichain $G_n$, $n\in \omega$. 
According to this theorem for every subset $S\subset \omega$ 
we must find an $\omega$-categorical 2-step nilpotent QE-group 
$G_S$ of exponent four such that $G_S$ is computable in $S$, 
and $G_d$ is embeddable into $G_S$ if and only if $d\in S$. 
For this purpose take the class $\mathcal{K}_S$ of all 
groups from $\mathcal{K}$ which do not embed all $G_d$ 
with $d\not\in S$.  
One easily sees that $\mathcal{K}_S$ is computable in $S$. 
On the other hand it is obvious that subgroups of groups 
from $\mathcal{K}_S$ belong to $\mathcal{K}_S$, and the 
free amalgamation defined for $\mathcal{K}$ guarantees 
the amalgamation (and the joint embedding) property for 
$\mathcal{K}_S$. 
Let $G_{S}$ be the Fra\"{i}ss\'{e} limit of the class 
$\mathcal{K}_S$. 
Consider axioms of $Th(G_S )$. 
As we already know we must formalize the following properties: 
\begin{quote}
(a) $\mathcal{K}_S$ coincides with the class of all 
finite substructures which can be embedded into $G_S$ and \\
(b) Every isomorphism between finite substructures of $G_S$ 
extends to an automorphism of $G_S$.   
\end{quote}
The first one is obviously formalized by $\forall$- and 
$\exists$-formulas and the set of these formulas is 
computable with respect to $S$. 
It is well-known that to formalize (b) we should express 
that for any two groups $H_1 <H_2$ from $\mathcal{K}_S$ 
any embedding of $H_1$ into $G_S$ extends to an embedding 
of $H_2$ into $G_S$. 
These sentences are $\forall\exists$ and obviously 
form a set computable in $S$ (in fact we may additionally 
assume that $H_2$ is 1-generated over $H_1$). 
As a result the theory $Th(G_S )$ is decidable in $S$. 
Thus it has a model computable in $S$. 
Since the theory is $\omega$-categorical we may assume 
that $G_S$ is computable in $S$. 

(2) We apply Theorem 2.3 from \cite{richter} 
to the effective antichain $G_n$, $n\in \omega$. 
According to this theorem for every subset $S\subset \omega$ 
we must find an $\omega$-categorical 2-step nilpotent QE-group 
$G_S$ of exponent four such that $G_S$ is enumeration reducible 
to $S$, and $G_d$ is embeddable into $G_S$ if and only if $d\in S$. 
For this purpose take the class $\mathcal{K}_S$ of all 
groups from $\mathcal{K}$ which do not embed all $G_d$ 
with $d\not\in S$ and repeat the construction of $G_S$ above.  

We now must additionally check that there is 
an effective procedure whose outputs enumerate 
$G_{S}$ when any enumeration of $S$ is supplied for the inputs.  
At the $n$-th step of an enumeration of $S$ we have a sequence 
$S_n =\{ s_0 ,...,s_n \}\subset S$. 
If $Q\subset G_S$ is the already enumerated part of $G_S$ 
let us consider all 1-types of $Th(G_S )$ over $Q$. 
By quantifier elimination they are quantifier free and 
the number of them depends on the isomorphism type of $Q$.   
At this step we choose (in turn) realizations of those  
types so that the subgroup generated by them together with $Q$ 
can be embedded into $G_{S_n}$. 
Since $S_n$ is finite, $Th(G_{S_n})$ is decidable. 
Thus this step can be done effectively. 

As a result we will obtain an enumeration of an elementary 
substructure of $G_S$. 
By model completeness and $\omega$-categoricity 
we see that it can be treated as an enumeration of $G_S$. 
$\Box$

\bigskip 

INSTITUTE OF MATHEMATICS, UNIVERSITY OF WROC{\L}AW, \parskip0pt 

pl.GRUNWALDZKI 2/4, 50-384 WROC{\L}AW, POLAND \parskip0pt

E-mail: 

ivanov@math.uni.wroc.pl 

biwanow@math.uni.wroc.pl

\end{document}